\documentclass[11pt]{article}
\usepackage{amssymb,amsmath,amsthm,amsfonts,verbatim,tikz,enumitem}
\usepackage[all,cmtip]{xy}
\usepackage[margin=1.25in]{geometry}
\usepackage{todonotes}
\usepackage{qtree}

\newtheorem{theorem}{Theorem}[section]

\newtheorem{conjecture}[theorem]{Conjecture}
\newtheorem{guiding}[theorem]{Guiding Principle}

\newtheorem{problem}[theorem]{Problem}
\newtheorem{question}[theorem]{Question}

\theoremstyle{definition}

\newtheorem{example}[theorem]{Example}
\newtheorem{examples}[theorem]{Examples}

\newtheorem{remark}[theorem]{Remark}
\newtheorem{remarks}[theorem]{Remarks}

\numberwithin{equation}{section}

\newcommand{\Poly}{{\rm Poly}}

\newcommand{\Ac}{\mathcal{A}}

\newcommand{\Mc}{\mathcal{M}}

\newcommand{\xdownarrow}[1]{%
  {\left\downarrow\vbox to #1{}\right.\kern-\nulldelimiterspace}
}

\newcommand{\Cb}{\mathbb{C}}

\newcommand{\Gb}{\mathbb{G}}

\newcommand{\Pb}{\mathbb{P}}

\newcommand{\Rb}{\mathbb{R}}

\newcommand{\Zb}{\mathbb{Z}}

\newcommand{\Mod}{Mod}

\DeclareMathOperator{\PConf}{PConf}

\DeclareMathOperator{\UConf}{UConf}

\newcommand{\para}[1]{\medskip\noindent\textbf{#1.}}

\DeclareMathOperator{\M}{{\cal M}}
\DeclareMathOperator{\N}{{\cal N}}

\newcommand\dual{\raise0.9ex\hbox{$\scriptscriptstyle\vee$}}

\title{Rigidity of moduli spaces and algebro-geometric constructions\\
{\em \normalsize In memory of S.S. Chern} \footnote{This paper is based on a lecture I gave at Chern's 110th birthday conference on October 12, 2021. I would like to thank Professors 
Shing-Tung Yau and Shiu-Yuen Cheng for soliciting the talk and this paper.  I also thank Dan Margalit for useful discussions.}}

\author{Benson Farb}

\begin{document}
\maketitle
\begin{abstract}
In this paper we propose two guiding principles that suggest a number of  
conjectures (some now proved) about various forms of rigidity for moduli spaces arising in algebraic geometry.  Such conjectures have group-theoretic, topological and holomorphic aspects, and so they also provide motivation for natural problems in geometric group theory and topology. 
\end{abstract}




\section{Introduction}

The purpose of this paper is to propose two guiding principles that suggest a number of  
conjectures (some now proved) about various forms of rigidity for moduli spaces arising in algebraic geometry.  Such conjectures have group-theoretic, topological and holomorphic aspects, and so they also provide motivation for natural problems in geometric group theory and topology.  These conjectures are in the style of, and in some cases can be seen as algebro-geometric interpretations of, rigidity theorems such as Margulis superrigidity. 

All varieties in this paper will be taken over $\Cb$, although the story over 
arbitrary fields would be quite interesting to explore.

\para{Moduli spaces from constructions}
One of the most compelling aspects of algebraic geometry is its remarkable constructions.   
Such constructions can often be reinterpreted as surprising morphisms between 
moduli spaces.  To choose three of many examples (see below for more): 

\begin{enumerate}
\item (Resolving the quartic) For each $n\geq 1$ let ${\rm Poly}_n$ be the space of monic, square-free, degree $n$ polynomials $P\in \Cb[x]$.  The space ${\rm Poly}_n$ is embedded 
$\Cb^n$ (by recording the coefficients of $P$ in $\Cb^n$) as the complement of a 
hypersurface, namely the zero set $Z(\Delta_n)$ of 
the discriminant polynomial $\Delta_n\in\Zb[a_1,\ldots,a_n]$.  The classical construction of ``resolving the quartic'' 
(Ferrari, 1545) produces a surjective morphism of quasiprojective varieties
\begin{equation}
\label{eq:quartic:resolve1}
{\rm R}: {\rm Poly}_4\to {\rm Poly}_3
\end{equation}
defined by taking the quartic polynomial $f\in{\rm Poly}_4$ with set of (distinct) roots $\{a_1,a_2,a_3,a_4\}$ to the cubic polynomial 
${\rm R}(f)\in\Poly_3$ with distinct \footnote{The miracle here is that if the $\{a_i\}$ are distinct then so are the $\{b_j\}$.} roots
%

\[b_1:=\frac{\displaystyle (a_1-a_2-a_3+a_4)^2}{\displaystyle 4}, \ \ 
b_2:=\frac{\displaystyle (a_1-a_2+a_3-a_4)^2}{\displaystyle 4}, \ \ 
b_3:=\frac{\displaystyle (a_1+a_2-a_3-a_4)^2}{\displaystyle 4}
\]

The induced map ${\rm R}_*:\pi_1({\rm Poly}_4)\to \pi_1({\rm Poly_3})$ induces a surjection 
of braid groups 
\[{\rm R}_*:B_4\to B_3.\]  By factoring through the surjection $B_n\to S_n$ to the symmetric group $S_n$, the map ${\rm R}$ provides an explanation for the existence of the exceptional surjection $S_4\to S_3$.  

\item (Jacobi varieties) The {\em Jacobian} construction produces from any genus $g\geq 1$ Riemann surface $X$ a $g$-dimensional principally polarized abelian variety
\[{\rm Jac}(X):=\Omega^1(X)^*/H_1(X;\Zb),\]
where $\Omega^1(X)^*$ is the dual of the space of holomorphic $1$-forms on $X$ and the embedding $H_1(X;\Zb)\to\Omega^1(X)^*$ is given by integration along a $1$-cycle.  This construction globalizes to the {\em period mapping}
\[{\rm Jac}:{\mathcal M}_g\to{\mathcal A}_g,\]
where ${\mathcal M}_g$ is the moduli space of genus $g$ Riemann surfaces and 
${\mathcal A}_g$ is the moduli space of principally polarized $g$-dimensional abelian varieties.

\item (Flex points and torsion) 
Fix $d,n\geq 1$.  Is there a way to choose in an algebraically (or even continuously) varying fashion an unordered set of $n$ distinct points on every smooth, degree $d$ plane curve?  
More precisely, let ${\mathcal C}_d$ be the parameter space of smooth, degree $d$ curves in $\Cb\Pb^2$, and let 
\[{\mathcal E}_{d,n}:=\{(C,\{z_1,\ldots,z_n\}): C\in {\mathcal C}_d, z_i\in C, z_i\neq z_j \ \text{if}\ i\neq j\}.\]
The projection $\pi_{d,n}(C,\{z_1,\ldots,z_n\}):=C$ induces a fiber bundle
\begin{equation}
\label{eq:bundle5}
\begin{array}{ccl}
{\rm UConf}_n(C)&\longrightarrow &{\mathcal E}_{d,n}\\
&&\Big\downarrow \pi_{d,n}\\
&&{\mathcal C}_d
\end{array}
\end{equation}
where ${\rm UConf}_n(C)$ denotes the space of configurations of unordered $n$-tuples of distinct points on $C$.  We are asking for an algebraic (resp.\ continuous) section of $\pi_{d,n}$.  
Remarkably, when $d=3$ and $n=9k^2, k\geq 1$, such algebraic sections exist, as we now explain.

Any (complex) smooth cubic curve $C\subset \Pb^2$ has exactly nine flex points ${\cal F}_C$.  Each $z\in {\cal F}(C)$ determines a unique abelian group structure on $C$ with identity $z$.  For each $k\geq 1$ the union ${\cal F}_{9k^2}(C)$ of these $k$-torsion points under these group structures produces an algebraic section of $\pi_{d,n}:{\mathcal E}_{d,n}\to {\cal C}_d$ 
defined by 
\[C\mapsto {\cal F}_{9k^2}(C).\]  Thus for each $k\geq 1$, any smooth cubic polynomial $f\in \Cb[x_1,x_2,x_3]$ comes equipped with a set of $9k^2$ solutions among its uncountably many solutions (here $C$ is the zero-set of $f$)!  
\end{enumerate}

The principle I want to present here predicts that each of the three maps constructed above 
is indeed special; each is unique in various senses.  For example, resolving the quartic is the only nontrivial way of producing, from a monic, square-free polynomial, another of lower degree, in a holomorphically (or even continuously) varying way (see \S\ref{section:polynomials} for details).   Another example is the following result, which I proved in \cite{Fa}.

\begin{theorem}[{\bf Global rigidity of the period mapping}]
\label{theorem:farb2}
Let $g\geq 3$ and assume that $h\leq g$.  Let $F:\Mc_{g}\to\Ac_h$ be any nonconstant holomorphic map of complex orbifolds.  Then $h=g$ and $F={\rm J}$.
\end{theorem}

See \S\ref{section:periods} for a discussion of this and related results and problems.

\para{Two guiding principles}  I'd like to try to formalize the statement that maps between moduli spaces (and sections of certain bundles) arising from algebro-geometric constructions are unique.  Such constructions typically start with a variety of a certain type, perhaps equipped with extra data such as a tuple of subvarieties, and then one produces a new variety, also perhaps with extra data, via geometry.  Such a construction induces a morphism, or more generally a rational map
\[F:\M\to\N\]
where $\M$ and $\N$ are the moduli spaces of the input/output varieties (perhaps equipped with their extra data) of the construction.  We call such a map between moduli spaces a {\em constructive map}.  While this is not a precise mathematical definition, to paraphrase supreme court justice Potter Stewart: 
I know a constructive map when I see it.\footnote{Stewart's comment described his threshold test for obscenity: ``I know it when I see it.''} We also apply this term to sections of bundles such as \eqref{eq:bundle5}, constructed via geometry or arithmetic.

Constructive maps are things of beauty.  Their existence often seems miraculous.  The following  is a way to formalize how special such constructions are. 

\begin{guiding}[{\bf Constructive maps are rigid}]
\label{guiding}
Constructive maps $F:\M\to\N$ satisfy the following: 
\begin{enumerate}
\item {\rm Group-theoretic rigidity: }Any nontrivial representation $\pi_1^{\rm orb}(\M)\to\pi_1^{\rm orb}(\N)$ is conjugate to $F_*$.
\item {\rm Topological rigidity: }any homotopically nontrivial continuous map $f:\M\to\N$ is homotopic to $F$.
\item {\rm Holomorphic rigidity: }any nonconstant holomorphic map $f:\M\to\N$ is equal \footnote{Sometimes we want to allow for pre-composition (resp.\ post-composition) with automorphisms of the domain (resp.\ codomain), for example by a linear change of coordinates.} to $F$. 
\item {\rm Rational rigidity: }any nonconstant rational map $f:M\to \N$ is equivalent as a rational map \footnote{Rational maps are {\em equivalent} as rational maps if they are equal on a Zariski open subset of the domain. Also see Footnote 3.}  to $F$.
\end{enumerate}
\end{guiding}

\begin{remarks}
\ 

\begin{enumerate}
\item If $\N$ is aspherical then topological rigidity and group-theoretic rigidity are equivalent.
\item One could replace ``holomorphic rigidity'' with ``morphism rigidity''; these are often equivalent (by extending holomorphic maps to certain compactifications and applying Chow's Theorem).
\item There are also natural variations and generalizations of Conjecture~\ref{guiding}, such as versions for maps between all finite covers of $\M$ and $\N$, and for certain special infinite covers (e.g.\ Torelli space).  
\end{enumerate}
\end{remarks}

The following is a natural generalization of Guiding Principle~\ref{guiding}.

\begin{guiding}[{\bf Characterization principle}]
Let $\{\M_i\}_{i\in I}$ and $\{\N_j\}_{j\in J}$ be natural families of moduli spaces.  Then any nonconstant holomorphic map (resp.\ morphism, rational map, etc.) $M_k\to\N_\ell$ with $k\in I, \ell\in J$ is constructive.
\end{guiding}

In terms of the examples above, the guiding principles predict:
\begin{enumerate}
\item Suppose $n>m>2$.  Let $F:\Poly_n\to\Poly_m$ be a nonconstant holomorphic map not factoring through the discriminant map $\Delta_n:\Poly_n\to\Cb^*$.  Then $(n,m)=(4,3)$ and $F$ is (up to affine coordinate changes in the domain and range) resolving the quartic.

\item Theorem~\ref{theorem:farb2}: 
Let $g\geq 3$ and assume that $h\leq g$.  Let $F:\Mc_{g}\to\Ac_h$ be any nonconstant holomorphic map of complex orbifolds.  Then $h=g$ and $F={\rm J}$.

\item For $(d,n)=(3,2)$, any algebraic section of \eqref{eq:bundle5} comes from a torsion construction; namely, is a sum of (variations of) the construction given - see  \S\ref{section:hypersurfaces} for the precise definition.
\end{enumerate}

I'd conjectured each of these results based on the guiding principles, and each was subsequently proved to be true; see \S\ref{section:polynomials}, \S\ref{section:periods} and \S\ref{section:hypersurfaces}.   A general reason in support of the holomorphic versions 
of the guiding principles is that holomorphic maps tend - with some notable exceptions - to be unique in their homotopy class.

I view these types of results as being worthwhile goals, for several reasons.  For one, they are in the same style as Mostow Rigidity, Margulis superrigidity, and other rigidity results of those 
types: characterizing a mathematical object within a larger class of such.  In proving such theorems, one is forced to understand the object in a much deeper way.  Secondly, I view attempts to prove these types of rigidity theorems as a kind of systematic search for beautiful algebro-geometric constructions.  Finally, trying to prove these conjectures brings up a host of interesting problems in group theory and topology, as we attempt to explain in this paper.

\para{General proof method}
In many cases there is a common outline to prove that a constructive map $F:{\mathcal M}\to{\mathcal N}$ is the unique 
holomorphic map ${\mathcal M}\to{\mathcal N}$. In the case when $\N$ is aspherical, the outline reads as follows: Let $f:{\mathcal M}\to{\mathcal N}$ be nonconstant and holomorphic.
\bigskip

\noindent
{\bf Step 1.} Prove group-theoretic rigidity. (Ideas on this later.)

\smallskip
\noindent
{\bf Step 2.} When $\N$ is aspherical this already implies that $f$ is either homotopically trivial or $f$ is topologically rigid.

\smallskip
\noindent
{\bf Step 3.} Now assume that $f$ is holomorphic.  If $\N$ is aspherical, then when it is homotopically trivial it is often constant by the maximum principal.  When $f$ is homotopically nontrivial, one can try to find enough ``negative curvature'' in $\N$ to prove uniqueness of holomorphic maps in a homotopy class. One example where one does {\em not} have actual Riemannian negative curvature but where such uniqueness holds is when $\N$ is the moduli space $\Mc_g$ of smooth, genus $g$ projective curves (see Theorem \ref{theorem:farb1} and the discussion following it).

\begin{remark}[{\bf Orbifold issues}]
Many of the moduli spaces discussed in this paper are (good) complex orbifolds: they are quotients $X/\Gamma$ of simply-connected manifolds $X$ by groups $\Gamma$ acting properly discontinuously but not necessarily freely.  A map $X/\Gamma\to Y/\Lambda$ between complex orbifolds is just an equivariant map $X\to Y$ with respect to a homomorphism $\rho:\Gamma\to\Lambda$.  When one of the group actions, say $\Gamma$ on $X$, is not faithful, then one can choose to replace $\Gamma$ by $\Gamma':=\Gamma/{\rm ker}(\rho)$, and 
ask for equivariance with respect to the $\Gamma'$ action.  This seemingly technical point can actually change the results; see \cite{Ser}.  
\end{remark}

\noindent
{\bf Outline of paper. }In the rest of this paper we give a sampling of problems, conjectures  and results motivated by the guiding principles. We restrict ourselves to four types of problems/conjectures: 
problems about polynomials (\S\ref{section:polynomials}); problems about choosing points on hypersurfaces (\S\ref{section:hypersurfaces}); problems 
about period mappings (\S\ref{section:periods}); and problems about constructions in enumerative geometry (\S\ref{section:enumerative}).   The choice of problems is of course biased, and is meant only to illustrate the guiding principles.  My hope is to inspire the reader to come up with their own rigidity conjectures, and hopefully even some theorems.

\section{Spaces of polynomials}
\label{section:polynomials}

Configuration spaces and maps between them provide a rich collection of 
examples of interesting maps to which to apply the guiding principles.   To set notation, given a topological space $X$ and an integer $n\geq 0$, let 
\[\PConf_n(X):=\{(x_1,\ldots,x_n) \in X^n: x_i\neq x_j \ \text{for}\ i\neq j\}\]
be the space of configurations of ordered $n$-tuples of distinct points in $X$.   The symmetric group $S_n$ acts on $X^n$ by permuting coordinates, leaving invariant 
the subspace $\PConf_n(X)$.  The quotient 
\[\UConf_n(X):=\PConf_n(X)/S_n\]
is the space of configurations of {\it unordered} $n$-tuples of distinct points in $X$.  Note that when $X$ is an algebraic variety then so are $\PConf_n(X)$ and $\UConf_n(X)$.  

The variety $\UConf_n(\Cb)$ relates to polynomials as follows.  For each $n\geq 1$, let ${\rm Poly}_n$ be the space of monic, square-free (i.e.\ has no repeated root), degree $n$ polynomials $P\in \Cb[x]$.   The classical theory of discriminants produces a polynomial $\Delta_n\in \Zb[x_1,\ldots,x_n]$ 
with the property that the polynomial $P(Z)=Z^n+a_1Z^{n-1}+\cdots +a_n$ is square-free if and only if $\Delta_n(a_1,\ldots,a_n)=0$.  It follows that 
\[{\rm Poly}_n=\{(a_1,\ldots,a_n): \Delta_n(a_1,\ldots ,a_n)\neq 0\}=\Cb^n -\{\Delta_n = 0\}\]
is a hypersurface complement in $\Cb^n$. We remark that the hypersurface $\{\Delta_n=0\}$ is highly singular and incredibly complicated: while $\Delta_2(b,c)=b^2-4c$, 
\[\Delta_3(b,c,d)=b^2c^2-4c^3-4b^3d-27d^2+18abc,\] and the complexity increases dramatically with $n$.  The Fundamental Theorem of Algebra implies that for each $n\geq 1$, the map 
\[{\rm Poly}_n\to \UConf_n(\Cb)\]
given by
\[P(Z)\mapsto \{\lambda: P(\lambda)=0\}\]
is an isomorphism of affine varieties.

The classical approach to understanding polynomials, initiated by Tschirnhaus in the 17th century, was to try to reduce formulas for the roots of the general degree $n$ polynomial to such formulas for polynomials of lower degrees. The most well-known example is resolving the quartic, which (as discussed above - see \eqref{eq:quartic:resolve1}) produces the morphism 
\begin{equation}
\label{eq:resolve6}
{\rm R}: {\rm Poly}_4\to {\rm Poly}_3.
\end{equation}  It is natural to show that no other such miracles were found because they do not actually exist.  Now, if one considers maps between spaces of all degree $n$ polynomials and degree $m$ polynomials, then each space is affine, and many such algebraic maps exist. Hence it is natural to restrict to the spaces ${\rm Poly}_n$ of square-free polynomials.  

Now, the ``resolving the quartic'' map $R$ is not quite unique.  First, for each $n\geq 1$ the discriminant gives a holomorphic surjection 
\[\Delta_n:{\rm Poly}_n\to \Cb^*,\]
and for any $m\geq 1$ we can compose $\Delta_n$ with any holomorphic map $\Cb^*\to {\rm Poly}_m$.  For $(n,m)=(4,3)$ we can also alter the  map $R$ of \eqref{eq:resolve6} 
by multiplying by a power $d\geq 0$ of the discriminant, giving a morphism 
\[{\rm R}_d: {\rm Poly}_4\to {\rm Poly}_3\]
defined by 
\[R_d(a_1,\ldots,a_n):=\Delta_n(a_1,\ldots,a_n)^dR(a_1,\ldots,a_n).\]
In particular $R_0=R$.  Note that if one replaced $\Delta_n(a_1,\ldots,a_n)^d$ by an arbitrary degree $d$ polynomial in $\Delta_n(a_1,\ldots,a_n)$ then this polynomial would take the value $0$ at some nonzero input, by the Fundamental Theorem of Algebra, and so the corresponding map $R_d$ would not be well-defined.   Up to these minor modifications, I conjecture that resolving the quartic is unique, in the following sense.

\begin{conjecture}[{\bf Resolving the quartic is unique}]
\label{conjecture:polymaps}
Let $m,n\geq 3$, and let $\Psi:{\rm Poly}_n\to {\rm Poly}_m$ be an arbitrary nonconstant holomorphic map. Assume that $\Psi$ is not constant, and that $\Psi$ does not factor through the discriminant $\Delta_n:{\rm Poly}_n\to \Cb^*$.  If $n\geq 4$ then $(n,m)=(4,3)$ and, after precomposing and postcomposing automorphisms of ${\rm Poly}_4$ and ${\rm Poly}_3$, respectively \footnote{These are induced by an affine map of $\Cb$, corresponding to a linear change of variables.}, it must be that $\Psi=R_d$ for some $d\geq 0$.
\end{conjecture}

The systematic study of holomorphic maps ${\rm Poly}_n\to {\rm Poly}_m$ was initiated by V. Lin, over a series of papers (see \cite{Li} and the references contained therein).  Lin proved in Theorem 3 of \cite{Li} that Conjecture~\ref{conjecture:polymaps} is true in the special case where $\Psi$ is of the form $\Psi(\{z_i\})\cap\{z_i\}=\emptyset$ \footnote{Here $\{z_i\}_{i=1}^n$ denotes a point in 
$\UConf_n(\Cb)$, which we are identifying with $\Poly_n$.}; that is, $\Psi$ only ``adds points'' (what Lin calls ``disjoining'').    Recent work of Chen-Kordek-Margalit \cite{CKM} solves the ``topological piece'' of Conjecture \ref{conjecture:polymaps} for the cases $5\leq m\leq 2n$ and $(n,m)=(4,3)$; see below.

The assumption $n\geq 4$ of Conjecture~\ref{conjecture:polymaps} is necessary, as there are many holomorphic maps ${\rm Poly}_3\to {\rm Poly}_m$ for various $m$, as we now explain.  Given a configuration $\lambda:=\{\lambda_1,\lambda_2,\lambda_3: \lambda_i\in\Cb, \lambda_i\neq \lambda_j\}$, let $E_\lambda$ be the $2$-sheeted branched cover of $\widehat{\Cb}$ branched over $\{\lambda_1,\lambda_2,\lambda_3,\infty\}$, each with ramification index $2$.  Then $E_\lambda$ is a complex torus, and comes equipped with a special point, the point at infinity, which we can take as a basepoint, endowing $E_\lambda$ with the structure of an elliptic curve with this basepoint as identity element and 
$\{\lambda_1,\lambda_2,\lambda_3,\infty\}$ as the set of $2$-torsion points.  The map 
\[\Psi_k:{\rm Poly}_3\to {\rm Poly}_{k^2-1}\]
defined by
\[\Psi_k(\{\lambda_1,\lambda_2,\lambda_3\}):=\{\text{nonzero $k$-torsion points of $E_\lambda$}\}\]
is a nonconstant  holomorphic map.  Note that $\Psi_k$ is disjoining for $k>2$.  Trevor Hyde (personal communication) has made progress on classifying all holomorphic maps 
${\rm Poly}_3\to {\rm Poly}_m$ with $m\geq 3$.  

One can ask a more general question: for fixed $m,n\geq 3$, what are the continuously (equivalently smoothly) varying ways (up to homotopy) of assigning a degree $m$ monic, square-free polynomial to a degree $n$ monic, square-free polynomial.  The following question was asked explicitly by the author and D. Margalit.  

\begin{question}
\label{question:continuous:polymaps}
Classify all homotopy classes of maps ${\rm Poly}_n\to {\rm Poly}_{n+k}$ for all $n\geq 3, k\geq 0$.  Equivalently, classify all homomorphisms $B_n\to B_{n+k}$, where $B_m$ denotes the braid group on $m$ strands.
\end{question}

The second part of \ref{question:continuous:polymaps} is equivalent to the first since ${\rm Poly}_n$ is a $K(B_n,1)$ space, where $B_n$ is the braid group on $n$ strands.  In addition to the obvious inclusions $B_n\to B_{n+k}$ for $k\geq 1$, there are, at least for $k\geq n$, many interesting homomorphisms $B_n\to B_{n+k}$, for example via braid cabling constructions, via ``transvections'', and via twists of these. Chen-Kordek-Margalit have a precise conjecture classifying all homomorphisms $B_n\to B_{n+k}$ for $n\geq 4$ and any $k\geq 1$.  However, the corresponding smooth maps ${\rm Poly}_n\to {\rm Poly}_{n+k}$ don't seem to be homotopic to holomorphic maps, although this is still open.  What do these maps correspond to at the level of polynomials? The case $k\leq n$ of Margalit's conjecture was solved by L. Chen-Kordek-Margalit as Theorem 1.1 of \cite{CKM}.  A first question that seems within reach:  do the cabling homomorphisms $B_n\to B_{rn}, r\geq 2$ correspond to an algebraic procedure on polynomials? Proposition 7.6 of \cite{CKM} classifies all surjections $B_4\to B_3$; the answer is compatible with Conjecture \ref{conjecture:polymaps}, since multiplying $\Delta_n^d$ induces on $\pi_1$ what they call a ``transvection''.

One can show, using the Schwarz Lemma, that if $\Psi:{\rm Poly}_n\to {\rm Poly}_m$ is holomorphic then $\Psi_*:B_n\to B_m$ takes each ``multitwist'' in $B_n$ to a ``multitwist'' 
in $B_m$.  It would be interesting to find other group-theoretic restrictions on $\Psi_*$ for holomorphic $\Psi$, compared to smooth $\Psi$.  As a basic illustration of what one can do with this type of constraint, recall the isomorphism of varieties $\Poly_n\cong \UConf_n(\Cb)$, and consider for $n\geq 2$ the continuous map 
\[\Phi_n:\UConf_n(\Cb)\to\UConf_{n+1}(\Cb)\]
defined by
\[\Phi_n(\{z_i\}):=\{z_i\}\cup \{\sum|z_i|+1\}.\]
The map $\Phi_n$ is not holomorphic, and I had conjectured that it is not even homotopic to a holomorphic map, except in the case $n=3$ discussed above.   This was proved by Hyde.

\begin{theorem}[Hyde]
\label{theorem:Hyde}
There is a morphism $\psi:\UConf_n(\Cb)\to\UConf_{n+1}(\Cb)$ if and only $n=3$.
\end{theorem}

One might hope for a similar theorem with target $\UConf_{n+k}(\Cb)$ for all $k\geq 1$.

\para{Polynomials with extra data}  As a  more refined question, one can ask: what are the smoothly (resp.\ holomorphically) varying ways of assigning, to a monic, square-free, degree $n$ polynomial together with some collections of its roots, a monic, square-free, degree $m$ polynomial.  In mathematical terms, given any subgroup $H<S_n$, consider the intermediate cover 
\[\PConf_n(\Cb)\to X_H\to {\rm Poly}_n\]
corresponding to $H$.  So, for example, when $H=1$ then $X_H=\PConf_n$ and when $H=S_{n-1}$ then 
\[X_H=\{(P,\lambda):P\in{\rm Poly}_n, P(\lambda)=0\}.\]

\begin{problem}[{\bf Polynomials with extra data}]
\label{problem:polymaps:general}
Fix $m,n\geq 3$.  Let $A<S_n$ and $B<S_m$ be any subgroups, and let $X_A$ (resp.\ $X_B$) be the cover of ${\rm Poly}_n$ (resp.\ ${\rm Poly}_n$) corresponding to $A$ (resp.\ $B$).  Classify all nonconstant holomorphic maps and all smooth homotopy classes of maps $X_A\to X_B$. 
\end{problem}

Problem~\ref{problem:polymaps:general} is still open in the case $A=B=1$.  Another interesting case is when $m=n$ and $A=1$ and $B=S_n$, so that $X_A=\PConf_n(\Cb)$ and 
$X_B={\rm \Poly}_n$.  Here there is the classical {\em Viete map} $V_n:\PConf_n(\Cb)\to \Poly_n$ sending a set of $n$ distinct complex numbers to the unique monic, degree $n$ polynomial with those roots:

\[V_n(r_1,\ldots,r_n):=(-\sigma_1(r_1,\ldots,r_n),\sigma_2(r_1,\ldots,r_n),\ldots,(-1)^n\sigma_n(r_1,\ldots,r_n))\]
where $\sigma_i$ is the $i^{\rm th}$ elementary symmetric polynomial: 
\[\begin{array}{ll}
\sigma_1(r_1,\ldots,r_n)&:=\sum_{i=1}^n r_i\\
\sigma_2(r_1,\ldots,r_n)&:= \sum_{1\leq i<j\leq n}r_ir_j\\
\vdots &\\
\sigma_n(r_1,\ldots,r_n)&:=r_1r_2\cdots r_n.
\end{array}\]

It seems fundamental to ask if this way of attaching a monic, degree $n$, square-free polynomial to a set of $n$ distinct complex numbers is unique.

\begin{conjecture}[{\bf The Viete map is unique}]
\label{conjecture:Viete:unique}
Fix $n\geq 4$.  Let $\Psi:\PConf_n(\Cb)\to \Poly_n$ be any nonconstant holomorphic map.  Assume further that $\Psi_*$ does not factor through an abelian group.  Then, up to a linear changes of variables (i.e.\ affine automorphisms of the domain and range), $\Psi$ is 
the Viete map $V_n$.
\end{conjecture}

Note that there are many homotopy classes of continuous maps $\PConf_n(\Cb)\to \PConf_n(\Cb)$, equivalently homomorphisms $P_n\to P_n$, where $P_n:=\pi_1(\PConf_n(\Cb))$ be the pure braid group on $n$ strands.  It fits into an exact sequence 
\[1\to P_n\to B_n\to S_n\to 1.\]
 One example of such a homomorphism is the `` forget some strands'' 
homomorphism $P_n\to P_{n-k}$ followed by a standard inclusion $P_{n-k}\to P_n$.  While the first homomorphism is induced by a holomorphic map, the second homomorphism 
is almost surely not (cf.\ Theorem~\ref{theorem:Hyde}).  Note that $P_3\cong F_2\times \Zb$, which surjects to $F_2$, and so one can produce many homomorphism 
\[P_n\twoheadrightarrow P_3\twoheadrightarrow F_2\to P_n.\]
However, one can argue that none of these homomorphisms is induced by a holomorphic map.  Thus Conjecture~\ref{conjecture:Viete:unique} survives.

\section{Choosing points on hypersurfaces}
\label{section:hypersurfaces}

Fix $d,n,N\geq 1$.  Let ${\mathcal C}_{d,N}$ be the parameter space of smooth, degree $d$ hypersurfaces in $\Cb\Pb^N$, and let 
\[{\mathcal E}_{d,N,n}:=\{(C,\{z_1,\ldots,z_n\}): C\in {\mathcal C}_{d,N}, z_i\in C, z_i\neq z_j \ \text{if}\ i\neq j\}.\]
The projection $\pi_{d,n}(C,\{z_1,\ldots,z_n\}):=C$ induces a fiber bundle
\begin{equation}
\label{eq:multisection1}
\begin{array}{ccl}
{\rm UConf}_n(C)&\longrightarrow &{\mathcal E}_{d,N,n}\\
&&\Big\downarrow \pi_{d,N,n}\\
&&{\mathcal C}_{d,N}
\end{array}
\end{equation}
where ${\rm UConf}_n(C)$ denotes the space of configurations of unordered $n$-tuples of distinct points on $C$. Note that ${\mathcal E}_{d,N,1}\to {\mathcal C}_{d,N}$ is the universal smooth, degree $d$ hypersurface in $\Pb^N$; for simplicity we denote it by $E_{d,N}$.  
A section of the bundle \eqref{eq:multisection1} is called an {\em $n$-multisection} of ${\mathcal E}_{d,N}\to {\mathcal C}_{d,N}$. 

\begin{question}[{\bf The point-choosing problem}]
\label{question:hypersurface:general}
Fix $d,N,n\geq 1$.  Is it possible to choose in a holomorphically (resp.\ continuously, algerbaically, rationally) varying manner an unordered $n$-tuple of distinct points on every smooth, degree $d$ hypersurface in $\Pb^N$? 

In other words: does there exist a holomorphic (resp.\ continuous, algebraic, rational) $n$-multisection of  ${\mathcal E}_{d,N}\to {\mathcal C}_{d,N}$?  When such a section exists, is it unique (resp.\ unique up to homotopy)?
\end{question}

Even the case $n=1$ of Question~\ref{question:hypersurface:general} is open (as far as I know).  Here is a first challenge:

\begin{conjecture}[{\bf You can't choose a point}]
Fix $d\geq 3, N\geq 2$. There is no holomorphic, or even continuous, section of ${\mathcal E}_{d,N}\to {\mathcal C}_{d,N}$.
\end{conjecture}

One of the compelling aspects of Question~\ref{question:hypersurface:general} is that there are many examples when such sections exist.  These (the holomorphic ones, anyway) seem always to arise from beautiful algebro-geometric and arithmetic constructions.  In trying to answer Question~\ref{question:hypersurface:general}, we are both proving uniqueness of these constructions (following the guiding principle), as well as initiating a systematic search for 
other, not yet known, algebro-geometric and arithmetic constructions.

\subsection{The case of cubic curves}  

It is a classical fact that every smooth cubic curve $C\subset\Pb^2$ has precisely $9$ {\em inflection points}; that is, points on $C$ whose tangent line intersects $C$ with multiplicity $3$.  
Let ${\rm Flex}(C)\in \UConf_9(\Pb^2)$ denote the set of flex points of $C$.  The map 
\[\Psi_9:{\mathcal C}_{3,2}\to {\mathcal E}_{3,2}\]
defined by 
\[\Psi_9(C):={\rm Flex}(C)\]
gives a $9$-multisection of ${\mathcal E}_{3,2}\to {\mathcal C}_{3,2}$.  Since the set of flex-points of $C$ is given by the intersection of $C$ with its associated Hessian curve (defined by the vanishing of the determinant of the Hessian of the equation defining $C$), it follows that $\Psi_9$ is an algebraic $9$-multisection.  Now, any choice of a point $p\in {\rm Flex}(C)$ gives $C$ the structure of an abelian group (elliptic curve) with $p$ as identity element and the set ${\rm Flex}(C)$ as the set of $3$-torsion points.  Thus $\{3\text{-torsion points}\}$ of $C$ is well-defined, independent of the choice of point in ${\rm Flex}(C)$.  For each $k\geq 1$ there is an algebraic 
$9k$-multisection $\Psi_{9k^2}$ of ${\mathcal E}_{3,2}\to {\mathcal C}_{3,2}$ defined by 
\[\Psi_{9k^2}:=\{\text{$3k$-torsion points of $C$}\}.\]
There are also (see, e.g. \cite{Ch}) algebraic multisections given by sets of $3k$-torsion points that are not $3j$-torsion points for any $j<k$, and indeed these have geometric interpretations.  For example, there are always precisely $72$ points on $C$ where an irreducible cubic intersects $C$ at points of multiplicity $9$; this set of points is the set of $9$-torsion points of $C$ that are not $3$-torsion points. We can also take unions of these torsion constructions to produce new algebraic multisections.  Banerjee-Chen call any such multisection a {\em multisection from torsion}.  As Chen observed in Theorem 1 of \cite{Ch} (atttributed to Maclaurin, Cayley and Gattazzo):

\begin{theorem}[{\bf Canonical subsets of smooth cubic curves}]
The universal cubic plane curve ${\mathcal E}_{3,2}\to {\mathcal C}_{3,2}$ admits an 
algebraic $n$-multisection precisely when 
\begin{equation}
\label{eq:multi:form}
n=9\sum_{m\in I}J_2(m)
\end{equation} 
where 
$J_2(m)=m^2\prod_{\text{$p$ prime, $p|m$}}(1-p^{-2})$ and $I$ is a set of positive integers.  For example when $n = 9, 27, 36, 72, 81, 99, 108, \ldots$.
\end{theorem}

Thus any smooth cubic plane curve $C$ has a canonical, algebraically varying (in $C$) set of $108$ points!  I found this to be remarkable, and conjectured that no other such constructions were possible, in the following precise senses:
\begin{enumerate}
\item Any algebraic $n$-multisection of ${\mathcal E}_{3,2}\to {\mathcal C}_{3,2}$ must be a multisection from torsion.
\item Any continuous $n$-multisection of ${\mathcal E}_{3,2}\to {\mathcal C}_{3,2}$ must be homotopic to a multisection from torsion; in particular no continuous $n$-multisections exist unless $n$ is of the form given in \eqref{eq:multi:form}.
\end{enumerate}

C. McMullen gave a proof of the first conjecture (see the Appendix of \cite{BC}).  Chen \cite{Ch} and Banerjee-Chen \cite{BC} proved many cases of the second conjecture, although many cases are still open.  As a sample of their results we mention:  if a continuous $n$-multisection exists then $9|n$; and a continuous $18$-multisection does not exist; and any continuous $9$-multisection is homotopic to $\Psi_9$.  Surprisingly, they found counterexamples to
the second conjecture:

\begin{theorem}[{\bf Banerjee-Chen \cite{BC}, Theorem 1.8}]
For any $m\geq 4$, the universal smooth cubic plane curve has a connected \footnote{A multisection is {\em connected}  if the space of smooth cubic plane curves equipped with a point in the multisection is connected.}, smooth $n$-multisection, $n=18J_2(m)$, that is not homotopic 
 to any multisection from torsion. 
\end{theorem}

The smallest examples of such an $n$ are $n=216, 432, 864, 1296, 2160, \ldots$.  The construction of such $n$-multisections in \cite{BC} is not geometric.  Is there a natural geometric construction of these? 

\subsection{Rational multisections} 
In contrast to the above situation, there are many examples of interesting {\em rational} multisections of the universal smooth, degree $d$ hypersurface ${\mathcal E}_{d,N}\to {\mathcal C}_{d,N}$; that is, multisections that are rational maps.  The point is that they are only defined on a Zariski open in 
${\mathcal C}_{d,N}$. Here are some examples.

\begin{enumerate}
\item The universal smooth quartic curve ${\mathcal E}_{4,2}\to {\mathcal C}_{4,2}$ has a rational $56$-multisection. Proof: It is classical that every smooth quartic curve $C\subset \Pb^2$ has $28$ {\em bitangents} - lines $L\subset\Pb^2$ with $L$ tangent to $C$ at two points (and so each intersection having multiplicity $2$).  These bitangents are not always distinct, but are counted with multiplicity.  For example, 
the Fermat quartic $F:=Z(x^4+y^4+z^4)\subset\Pb^2$ has $12$ standard bitangents, as well as $16$ points where the tangent to $F$ intersects $F$ in a single point of multiplicity $4$, giving a total of $40$ tangency points.  

The condition for a smooth quartic to have a point of tangency of multiplicity $4$ is a polynomial condition, and so there is a Zariski open $U\subset {\mathcal C}_{4,2}$  for which each $C\in U$ has $28$ distinct bitangents, and so $56$ points of tangency.  This gives a rational 
$56$-multisection of ${\mathcal E}_{4,2}\to {\mathcal C}_{4,2}$ with $U$ as domain of definition.

\item A general smooth, degree $d\geq 3$ plane curve $C$ has $3d(d-2)$ flex points, counted with appropriate multiplicity; for $d>3$ one must pass to a Zariski open subset 
$U\subsetneq {\mathcal C}_{d,2}$ in order for each curve $C\in U$ to have $3d(d-2)$ distinct flex points.   This gives a rational $3d(d-2)$ multisection of ${\mathcal E}_{d,2}\to {\mathcal C}_{d,2}$.  There is a similar story for bitangents of smooth plane curves of degree $d>4$.  

\item Each smooth cubic surface $S\subset\Pb^3$ contains $27$ distinct lines.  There is a Zariski open subset $U$ of the parameter space ${\mathcal C}_{3,3}$ of smooth cubic surfaces whose collection of lines have $135$ points of intersection.  This gives a rational $135$-section of the universal smooth cubic surface ${\mathcal E}_{3,3}\to {\mathcal C}_{3,3}$.  This section is not a morphism since $U\subsetneq {\mathcal C}_{3,3}$, due to the existence of so-called {\em Eckhardt points}, where three of the $27$ intersect. 

\item Here is a more mundane example.  Fix $d,r\geq 1$ and $N\geq 2$.  For each smooth, 
degree $r$ curve $C\subset \Pb^N$, there is a Zariski open $U_C\subset {\mathcal C}_{d,N}$ of consisting of those degree $d$ smooth hypersurfaces in $\Pb^N$ intersecting $C$ in $dr$ points.  This gives a rational $dr$-multisection of ${\mathcal E}_{d,N}\to {\mathcal C}_{d,N}$.
\end{enumerate}

The following seems quite natural, although I do not know how difficult it is.

\begin{problem}[{\bf Classification problem for rational multisections}]
Classify all rational multisections of ${\mathcal E}_{d,N}\to {\mathcal C}_{d,N}$.
\end{problem}

\section{Period mappings}
\label{section:periods}

A common classical construction is to associate to a variety, which is fundamentally nonlinear, a variety constructed via linear algebra, such as an abelian variety. This is typically achieved via Hodge theory. In this section we start by discussing the more classical case of (generalized) Jacobians, and continue with more general period domains.  We will then indicate how the uniqueness results one might hope for by applying the guiding principles are closely related to rigidity results for lattices in semisimple Lie groups.

\subsection{Jacobians}
The most classical way to attach one variety to another in a canonical way is the Jacobian construction, going back at least to Riemann (if not Abel).

\begin{examples}[{\bf Jacobians}]
\label{example:jacobians1}
Many period mapping are Jacobians.  Here are some examples; see Debarre's survey \cite{De} for many more.
\begin{enumerate}
\item The classical period mapping ${\rm J}:{\cal M}_g\to{\cal A}_g$, as described in the introduction.  An equivalent formulation, using notation that will be more suggestive for what will soon come: if $X$ is a smooth curve with Hodge decomposition $H^1(X;\Cb)=H^{1,0}(X)\oplus  H^{0,1}(X)$, then 
\[{\rm J}(X):=H^{0,1}(X)/H^1(X;\Zb).\]
\item Let ${\mathcal C}_{3,4}$ be the parameter space of smooth cubic hypersurfaces in $\Pb^4$.  If $X\in{\mathcal C}_{3,4}$ then Hodge theory implies $H^{1,2}(X)\cong\Cb^5$, 
and \[{\rm J}_{3,4}(X):=H^{1,2}(X)/H^3(X;\Zb)\]
is a $5$-dimensional principally polarized abelian variety.  This gives a holomorphic map 
\[{\rm J}_{3,4}:{\mathcal C}_{3,4}\to {\mathcal A}_{5}\]
called the {\em intermediate Jacobian}.
\item Let $X\in{\mathcal M}_g$ be a smooth, genus $g$ curve.  Any nonzero $\theta\in H^1(X;\Zb/2\Zb)$ determines an unbranched double cover $p:Y\to X$ with deck transformation $\sigma$.  Note that $Y$ is a genus $2g-1$ curve.  The {\em Prym variety} ${\rm Prym}(X,\theta)$ associated to $(X,\theta)$ is defined to be 
\[{\rm Prym}(X,\theta):=\frac{\displaystyle J(X)}{\displaystyle p^*(J(X))}. \]
Note that ${\rm Prym}(X,\theta)\in {\mathcal A}_{g-1}$.  Let
\[{\mathcal R}_g:=\{(X,\theta): X\in{\mathcal M}_g \ 
\text{and}\ 0\neq\theta\in H^1(X;\Zb/2\Zb)\}/\sim\]
where $(X_1,\theta_1)\sim(X_2,\theta_2)$ if there is an isomorphism $f:X_1\to X_2$ such that $f^*\theta_2=\theta_1$.  Then ${\mathcal R}_g$ is a quasiprojective variety and the map
\[{\rm Prym}:{\mathcal R}_g\to {\mathcal A}_{g-1}\]
is a morphism.  
\end{enumerate}
\end{examples}

\bigskip
Such constructions inspire the following.

\begin{question}[{\bf Uniqueness of Jacobians}]
Given a moduli space ${\mathcal M}$ of varieties, is there a nontrivial way to attach in a holomorphically (or even continuously) varying way a principally polarized abelian variety to each $X\in {\mathcal M}$?  More precisely, does there exist $g\geq 1$ and a nonconstant holomorphic map (or nontrivial homotopy class of continuous maps) $F:{\mathcal M}\to {\mathcal A}_g$? Is each such homotopy class represented by a holomorphic map? If so, is such a map unique?
\end{question}

The guiding principles indicate that classical constructions as in 
Example~\ref{example:jacobians1} should be unique.     I proved in \cite{Fa} that this is the case for the classical Jacobian.

\begin{theorem}[{\bf Global rigidity of the period mapping}]
\label{theorem:farb1}
Let $g\geq 3$ and assume that $h\leq g$.  Let $F:\Mc_{g}\to\Ac_h$ be any nonconstant holomorphic map of complex orbifolds.  Then $h=g$ and $F={\rm J}$.
\end{theorem}

As with many of the results and problems discussed in the present paper, the statement of Theorem~\ref{theorem:farb1} really has two parts:

\begin{enumerate}
\item (Topological rigidity) If $g\geq 3$ and $h\leq g$ and if $F:\Mc_{g}\to\Ac_h$ is any continuous map of topological orbifolds, then either $F$ is homotopically trivial or $h=g$ and $F$ is homotopic to the classical period mapping ${\rm J}$.  Since ${\cal A}_g$ is the quotient of the lattice ${\rm Sp}(2g,\Zb)$ acting on a nonpositively curved (hence contractible!) symmetric space, topological rigidity reduces to understanding homomorphisms ${\rm Mod}(S_g)\to {\rm Sp}(2g,\Zb)$.  This was essentially done by Korkmaz.
\item (Holomorphic rigidity) The general philosophy is that holomorphic maps are typically unique in their homotopy class unless their images lie in a product.  When targets (such as ${\mathcal A}_g$) are finite volume quotients of bounded symmetric domains, Borel-Narasimhan reduce such uniqueness in a homotopy class to having a single image point in common.  Finding such a point can be involved and difficult.  This is accomplished in \cite{Fa} using what I called the ``Wirtinger squeeze'':  a convexity argument applied to a Wirtinger-type inequality proves that a homotopy of holomorphic maps must restrict to every curve in ${\mathcal M}_g$ as an {\em algebraic} deformation.  Using a criterion of Saito, a rigid curve $C$ is then proved to exist, providing the needed point (any point of $C$).
\end{enumerate}

I conjectured in \cite{Fa} that such a result should hold for the Prym construction. This conjecture was recently proven by C. Servan in \cite{Ser}

\begin{theorem}[Uniqueness of the Prym construction]
\label{theorem:servan}
Let $g\geq 4$ and let Let $1\leq h\leq g-1$.  Let $F:{\mathcal R}_g\to {\mathcal A}_h$ be any nonconstant holomorphic maps of complex orbifolds.  Then $h=g-1$ and $F={\rm Prym}$.
\end{theorem}

In order to prove the topological rigidity part of Theorem~\ref{theorem:servan}, Servan classifies all homomorphisms 
\[{\rm Stab}_{{\rm Mod}(S_g)}(\theta)\to {\rm Sp}(2h,\Zb),\]
where $0\neq \theta\in H^1(S_g;\Zb/2\Zb)$ and ${\rm Stab}_{{\rm Mod}(S_g)}(\theta)$ is the stabilizer of $\theta$ of the natural action of $\Mod(S_g)$ on $H^1(S_g;\Zb/2\Zb)$.  To do this he classifies all homomorphisms
\[{\rm Stab}_{{\rm Mod}(S_{2g-1})}(\sigma) \to {\rm Sp}(2h,\Zb).\]
This and many other such questions provide interesting problems in group theory and low-dimensional topology.

Extending all of the above, it is natural to ask for ways to attach one curve to another, or a principally polarized abelian variety (of arbitrary dimension) to a curve.  One obvious construction of this is the following: fix an unramified, characteristic cover $p:S_h\to S_g$ of surfaces.  Any complex structure on $S_g$ can be pulled back via $p$ to a complex structure on $S_h$, thus inducing a holomorphic map
\[p^*:\Mc_g\to\Mc_h\]
and hence a holomorphic map $\psi:\Mc_g\to\Ac_h$ given by $\psi:{\rm J}\circ p^*$.  We call such maps {\em covering constructions}.  Note that we needed to choose a characteristic cover so as not to make any choices, thus giving a map with domain $\Mc_g$, as opposed to (for example) the moduli space ${\mathcal R}_g$ of Pryms.  C. McMullen asked me the following:

\begin{question}[{\bf Curves from curves}]
\label{question:curt}
Do all nonconstant holomorphic maps $\Mc_g\to\Mc_h$ and $\Mc_g\to\Ac_h$ come from covering constructions?
\end{question}

The smallest characteristic cover of $S_g$ is probably (it would be good to check this) the mod $2$ homology cover of $S_g$; in this case $h=2^{2g}-1$.  In particular, if the answer to Question~\ref{question:curt} is ``yes'', any holomorphic map $\Mc_g\to\Mc_h$ with $h<2^{2g}-1$ should be constant.  The best known result in this direction is due to 
Antonakoudis-Aramayona-Souto \cite{AAS}, who proved this statement for $h\leq 2g-2$.  

One can do a similar construction, and ask a similar question as Question~\ref{question:curt}, for branched covers, in which case one must consider moduli spaces of pointed curves.

\subsection{Other period mappings and superrigidity}

Many period mappings take values in other finite volume quotients of bounded symmetric domains.

\begin{example}[{\bf The period mapping for smooth quartic surfaces}] 
\label{example:quartic4}
Let $X\subset\Pb^3$ be a smooth quartic surface.  The Hodge decomposition of the primitive cohomology of any such $X$ is 
\begin{equation}
\label{eq:quartic3}
H^2(X;\Cb)_{\rm prim}=H^{0,2}(X)\oplus H^{1,1}(X)_{\rm prim}\oplus H^{2,0}(X)\cong \Cb\oplus \Cb^{19}\oplus\Cb
\end{equation}
and the intersection form $Q_X$ on $H^2(X;\Cb)_{\rm prim}$ has signature $(2,19)$.  It follows that the Hodge decomposition \eqref{eq:quartic3} is determined by the line $H^{2,0}(X)$ in $H^2(X;\Cb)_{\rm prim}$; that is, by a point in the projectivization $\Pb(H^2(X;\Cb)_{\rm prim})\cong \Pb^{20}$.    The period mapping takes values in the quotient of the bounded symmetric domain 
\[\{[\omega]\in\Pb^{20}: Q(\omega,\omega)=0 \ \text{and}\ Q(\omega,\bar{\omega}>0\}
\cong {\rm SO}(2,19)^0/{\rm SO}(2)\times {\rm SO(19)}\]
by the cofinite volume arithmetic lattice ${\rm Aut}(H^2(X;\Cb)_{\rm prim}\cap H^2(X;\Zb), Q_X)$, which is a locally symmetric quasiprojective variety which we will denote by $M_{2,19}$.  The period map is thus a morphism
\[\Psi:  {\mathcal C}_{4,3}\to M_{2,19}\]
where ${\mathcal C}_{4,3}$ is the parameter space of smooth quartic surfaces $X\subset\Pb^3$.
\end{example}

More generally, the target of period mappings is complex manifold $N$ that is a fiber bundle $C\to N\to M$ where $C$ is a compact homogeneous space and $M$ is a finite volume locally symmetric (but not necessarily Hermitian!) locally symmetric variety.  Period mappings are constructed from Hodge-theoretic data, which is linear, preserving either a symmetric bilinear form on a lattice or a symplectic form (we have seen above examples of each). In either case the targets of period mappings are compact homogeneous space bundles over finite volume, locally symmetric varieties.   These are called {\em period domains} \footnote{Sometimes this term is reserved for their universal covers.}  

It is natural to ask if there are other, Hodge-theoretic or not, linear data to attach to families of varieties.

\begin{question}[{\bf Uniqueness of period mappings}]
\label{question:periods2}
Let ${\cal M}$ be a moduli space of smooth varieties.  Is the standard period mapping on 
${\cal M}$ unique? More precisely, what are the nonconstant holomorphic maps 
(resp.\ nontrivial homotopy classes of continuous maps) ${\cal M}\to N$ where $N$ is a period domain?
\end{question}

It seems to me that Question~\ref{question:periods2} is closely related to various rigidity theorems for lattices in real semisimple Lie groups; in particular it can be viewed as a generalization of Margulis's superrigidity theorem.  Superrigidity states that if $\Gamma$ is an irreducible lattice in 
a real semisimple Lie group $G$ with finite center, no compact factors and ${\rm rank}_{\Rb}(G)\geq 2$, then any linear representation of $G$ is either precompact or virtually extends to an algebraic representation of $G$.  Note that locally symmetric spaces are precisely products of compact symmetric spaces, flat tori, and space $\Gamma\backslash G/K$, where $G$ is a semisimple Lie group with no compact factors, $K$ is a maximal compact subgroup of $G$, and $\Gamma$ is a lattice in $G$ (that is, a cofinite volume discrete subgroup).  

I think it would be interesting and useful to interpret superrigidity as a kind of algebro-geometric rigidity. For example, the answer to the question: ``what are the ways of attaching one PPAV to another (perhaps of a different dimension)?'' can be worked out by applying superrigidity \footnote{Actually, for this case superrigidity is due to Bass-Milnor-Serre.} to homomorphisms ${\rm Sp}(2g,\Zb)\to {\rm Sp}(2h,\Zb)$? One can (and should) extend this to level structures as well.  

As another example, recall that the {\em Kummer construction} takes as input an abelian surface $A=\Cb^2/\Lambda \in\Ac_2$, blows up $A$ at its sixteen $2$-torsion points, and takes the quotient of the blowup by the involution induced by $(z,w)\mapsto (-z,-w)$, giving a K3 surface ${\rm Kum}(A)$.  The map $A\mapsto {\rm Kum}(A)$ induces a holomorphic map of moduli spaces 
\[{\rm Kum} : \Ac_2\to M_{2,19}\]
where $M_{2,19}$ is the locally symmetric variety defined above.   I believe that one can prove that the Kummer construction is unique among many possible constructions.  As a first step to proving this, one can check that superrigidity and some representation theory imply 
that any representation 
\[\rho: \pi_1^{\rm orb}(\Ac_2)\cong {\rm Sp}(4,\Zb)\to {\rm SO}(2,19)(\Rb)^0\]
is (conjugate to) the standard one, induced by ${\rm Kum}$.  

\section{Constructive maps via enumerative geometry}
\label{section:enumerative}

Various classical constructions from enumerative algebraic geometry can be viewed as highly nontrivial - sometimes surprising - morphisms of moduli spaces. One can then try to characterize these morphisms in various ways.  We state here only one example, and encourage the reader to find others.   

Recall the notation above that ${\mathcal C}_{d,N}$ is the parameter space of smooth, degree $d$ hypersurfaces in $\Pb^n$.  Denote by $\Gb(k,n)\cong {\rm Gr}(k+1,n+1)$ 
the Grassmannian of projective $k$-planes in $\Pb^n$.  

The Cayley-Salmon Theorem states that every smooth cubic surface $X\subset\Pb^3$ contains $27$ distinct lines.  The theory behind this implies that the map 
\begin{equation}
\label{eq:cubic:conf}
\Psi:{\mathcal C}_{3,3}\to \UConf_{27}(\Gb(1,3))
\end{equation}
is a morphism of quasiprojective varieties.  Similarly, every smooth quartic curve $C\subset\Pb^2$ has exactly $28$ {\em bitangents}; that is, lines in $\Pb^2$ tangent to $C$ with total multiplicity $4$ (i.e.\ in two points of multiplicity $2$ or one of multiplicity $4$).  
This gives a morphism of quasiprojective varieties
\begin{equation}
\label{eq:quartic:conf}
\Phi:{\mathcal C}_{4,2}\to \UConf_{28}(\Gb(1,2))
\end{equation}

These seem like a really interesting maps.  Since any two cubic surfaces in $\Pb^3$ intersect in a curve of degree at least $9$ by Bezout, it follows that the unordered set of $27$ lines determines the cubic surface; that is, that $\Psi$ is injective. A similar argument shows that $\Phi$ is injective.  The guiding principles suggest that there are no other ways to attach configurations of any other number of lines to every smooth cubic surface, or to quartic curves.  One might try to take all lines connecting all possible intersections of the $27$ lines, but this won't work since this number can vary over ${\mathcal C}_{3,3}$ (cf. the existence of Eckhardt points, mentioned above).  There is a similar problem for the bitangents.

\begin{conjecture}[{\bf Uniqueness Conjecture}]
Let $r,s\geq 1$.  Let $F:{\mathcal C}_{3,3}\to \UConf_{r}(\Gb(1,3))$ and 
$G:{\mathcal C}_{4,2}\to \UConf_{s}(\Gb(1,2))$ be any continuous maps.   
\begin{enumerate}
\item $F$ is homotopically trivial unless $r=27$ and $F$ is homotopic to the map $\Psi$ in \eqref{eq:cubic:conf}.  If $F$ is nonconstant holomorphic then $r=27$ and $F=\Psi$.
\item $G$ is homotopically trivial unless $s=28$ and $G$ is homotopic to the map $\Phi$ in \eqref{eq:quartic:conf}.  If $G$ is nonconstant holomorphic then $s=28$ and $G=\Phi$.
\end{enumerate}
\end{conjecture}

\bigskip{\noindent
Dept. of Mathematics\\ 
University of Chicago\\
E-mail: bensonfarb@gmail.com

\end{document}